\documentclass[12pt]{amsart}

\usepackage{latexsym,amsmath,amsthm,amscd,amsfonts,diagrams,graphics}

\DeclareFontFamily{U}{rsf}{}
\DeclareFontShape{U}{rsf}{m}{n}{
  <5> <6> rsfs5 <7> <8> <9> rsfs7 <10-> rsfs10}{}
\DeclareMathAlphabet{\mathscr}{U}{rsf}{m}{n}

\DeclareMathAlphabet{\mathgth}{U}{euf}{m}{n}

\DeclareFontFamily{U}{cyr}{}
\DeclareFontShape{U}{cyr}{m}{n}{
  <5> wncyr5 <6> wncyr6 <7> wncyr7 <8> wncyr8 <9> wncyr9 <10-> wncyr10}{}
\DeclareMathAlphabet{\mathcyr}{U}{cyr}{m}{n}

\input cyracc.def

\newcommand{\Sha}{\mathcyr{\cyracc Sh}}

\newcommand{\cA}{{\mathscr A}}

\newcommand{\cal}{\mathcal}

\newcommand{\sL}{{\cal L}}

\newcommand{\cE}{{\mathscr E}}
\newcommand{\cF}{{\mathscr F}}
\newcommand{\cG}{{\mathscr G}}

\newcommand{\cL}{{\mathscr L}}
\newcommand{\cM}{{\mathscr M}}

\newcommand{\cO}{{\mathscr O}}
\newcommand{\cP}{{\mathscr P}}
\newcommand{\cQ}{{\mathscr Q}}

\newcommand{\cX}{{\mathscr X}}

\newcommand{\balpha}{{\bar{\alpha}}}
\newcommand{\bbeta}{{\bar{\beta}}}

\newcommand{\tH}{\tilde{H}}

\newcommand{\dcE}{{\cE^{\cdot}}}

\newcommand{\FMYX}{\Phi_{Y\ra X}}

\newcommand{\fmXM}{\phi_{X\ra M}}
\newcommand{\FMXM}{\Phi_{X\ra M}}
\newcommand{\FMMX}{\Phi_{M\ra X}}

\newcommand{\sHom}{\underline{\mathrm{Hom}}}

\newcommand{\bProj}{\mathbf{Proj}}
\newcommand{\D}{{\mathbf D}_{\mathrm{coh}}^b}

\newcommand{\chk}{{\scriptscriptstyle\vee}}
\newcommand{\R}{\mathbf{R}}
\newcommand{\Ld}{\mathbf{L}}
\newcommand{\lotimes}{\stackrel{\Ld}{\otimes}}
\newcommand{\et}{{\mathrm{\acute{e}t}}}
\newcommand{\an}{{\mathrm{an}}}
\newcommand{\tors}{{\mathrm{tors}}}

\newcommand{\Kt}{{\mathrm{K3}}}

\DeclareMathOperator{\Hom}{Hom}

\DeclareMathOperator{\Ker}{Ker}

\DeclareMathOperator{\Pic}{Pic}

\DeclareMathOperator{\PGL}{PGL}
\DeclareMathOperator{\GL}{GL}

\DeclareMathOperator{\SL}{SL}
\DeclareMathOperator{\Td}{td}
\DeclareMathOperator{\ch}{ch}

\DeclareMathOperator{\id}{id}

\DeclareMathOperator{\Br}{Br}

\DeclareMathOperator{\NS}{NS}

\DeclareMathOperator{\rk}{rk}
\DeclareMathOperator{\Ext}{Ext}

\newcommand{\ra}{\rightarrow}
\newcommand{\lra}{\longrightarrow}

\newcommand{\scdot}{\,\cdot\,}

\newcommand{\C}{\mathbf{C}}

\newcommand{\Q}{\mathbf{Q}}
\newcommand{\Z}{\mathbf{Z}}
\newcommand{\gMod}{\mathgth{Mod}}
\newcommand{\gCoh}{\mathgth{Coh}}

\newcommand{\iso}{\cong}
\newcommand{\Cech}{\v{C}ech }

\newcommand{\pj}{\mathbf{P}}

\newarrow{Equal} =====

\theoremstyle{plain}
\newtheorem{theorem}{Theorem}[section]

\newtheorem{proposition}[theorem]{Proposition}
\newtheorem{conjecture}[theorem]{Conjecture}
\theoremstyle{definition}

\newtheorem{definition-theorem}[theorem]{Definition-Theorem}

\theoremstyle{remark}
\newtheorem{remark}[theorem]{Remark}

\renewcommand{\phi}{\varphi}

\begin{document}

\author{Andrei C\u ald\u araru}

\title{Non-Fine Moduli Spaces of Sheaves on K3 Surfaces}

\date{}

\begin{abstract}
In general, if $M$ is a moduli space of stable sheaves on $X$, there
is a unique $\alpha$ in the Brauer group of $M$ such that a $\pi_M^*
\alpha^{-1}$-twisted universal sheaf exists on $X\times M$.  In this
paper we study the situation when $X$ and $M$ are K3 surfaces, and we
identify $\alpha$ in terms of Mukai's map between the cohomology of
$X$ and of $M$ (defined by means of a quasi-universal sheaf).  We
prove that the derived category of sheaves on $X$ and the derived
category of $\alpha$-twisted sheaves on $M$ are equivalent.

This suggests a conjecture which describes, in terms of Hodge
isometries of lattices, when derived categories of twisted sheaves on
two K3 surfaces are equivalent.  If proven true, this conjecture would
generalize a theorem of Orlov and a recent result of Donagi-Pantev.
\end{abstract}

\maketitle

\section{Introduction}

\subsection{}
Ever since Mukai's seminal paper~\cite{Muk}, moduli spaces of sheaves
on K3 surfaces have attracted constant interest, largely because of
the wealth of geometric properties they possess.  Such moduli spaces
are always even-dimensional, so the first case with non-trivial
geometry is when $M$ is a 2-dimensional moduli space of sheaves on a
K3 surface $X$.  In this paper we'll focus on this situation, whose
study was initiated by Mukai, and later continued in a slightly
different direction by Orlov~(\cite{Orl}).

The guiding principle when studying such moduli spaces can be loosely
stated by saying that the geometry of $M$ is inherited from $X$.  In
our case, for example, $M$ is again a K3; however, $M$ could be quite
different from $X$, and this is why this situation is so interesting.
The question that arises is what relations can be found between
geometric objects on $X$ and on $M$, and most of the time we will be
interested in relating the cohomology groups (integral or rational,
plus Hodge structure) of $X$ and of $M$, as well as certain derived
categories of sheaves on $X$ and $M$.

\subsection{}
Mukai's original paper was concerned with studying the cohomology of
$M$, and we'll start by presenting a sketch of his ideas.  This can be
done most easily in the case when $M$ is {\em fine}, in other words
when a universal sheaf exists on $X\times M$.

The {\em Mukai lattice} of $X$ is defined as the group
\[ \tH(X,\Z) = H^0(X, \Z) \oplus H^2(X, \Z) \oplus H^4(X, \Z), \]
endowed with the Mukai product~(\ref{subsec:MukProd}), which is a
slight modification of the usual product in cohomology.  Using the
Chern classes of a universal sheaf one defines~(\ref{subsec:MukIsom})
a class in $H^*(X\times M, \Z)$, which induces a correspondence
\[ \phi:\tH(X, \Z) \ra \tH(M, \Z). \]
This correspondence is in fact an isomorphism, which respects the
Mukai product and the natural Hodge structure on the Mukai lattice
induced from the Hodge structure of $X$~(\ref{subsec:MukHodge}).  Thus
the integral Mukai lattices of $X$ and of $M$ are Hodge isometric, and
from this Mukai is able to give a complete description of $M$, in
terms of the Torelli theorem.

\subsection{}
Although it can be presented entirely in terms of cohomology, the
above construction relies heavily on the fact that there is an
underlying equivalence of derived categories between $X$ and $M$.
More precisely, the integral functor defined using a universal sheaf
is a Fourier-Mukai transform, i.e.\ an equivalence $\D(M)\iso \D(X)$.
It is an interesting question to ask for what other K3's $M$ do we
have such equivalences, and Orlov's work has provided an answer to
this question.  His main result can be stated as follows:
\vspace{2mm}

\noindent
{\bf Theorem (Orlov~\cite{Orl}).}  
{\em 
Let $X$ and $M$ be K3 surfaces.  Then the following are equivalent:
\begin{enumerate}
\item
$M$ is a fine, compact, 2-dimensional moduli space of stable sheaves
on $X$;
\item 
there is a Hodge isometry $T_X\iso T_M$ between the transcendental
lattices of $X$ and of $M$;
\item
the derived categories $\D(X)$ and $\D(M)$ are equivalent.
\end{enumerate} }

The implication $(1)\Rightarrow(2)$ is an immediate consequence of
Mukai's results: since the map $\phi$ is defined by means of an
algebraic class, it takes algebraic classes on $X$ to algebraic
classes on $M$.  Being a Hodge isometry, it follows that the
transcendental lattice of $X$ will be mapped by $\phi$ Hodge
isometrically onto the transcendental lattice of $M$.

\subsection{}
\label{subsec:sublatcat}
Orlov's theorem should be viewed as a derived version of the Torelli
theorem: given a Mukai lattice with a Hodge structure, singling out
the $H^2$ lattice will determine the K3 surface (by Torelli), and
therefore the category of coherent sheaves; singling out a smaller
lattice as the transcendental lattice of a K3 will no longer determine
the surface itself, but it will determine the derived category of
sheaves on the K3.

\subsection{}
The issue with the above results is the fact that the condition that
$M$ be fine is quite restrictive (at the moment, only one explicit
class of examples of this type is known).  There are many situations
in which relaxing this condition would be relevant -- for example, the
first example studied by Mukai~(\cite[2.2]{MukK3}) is not fine.  (In
this case one calculates the moduli space of spinor bundles on a
$(2,2,2)$ complete intersection in $\pj^5$ and finds it to be a double
cover of $\pj^2$ branched over a sextic.)  On the level of cohomology,
Mukai was able to avoid the condition of fineness, defining the map
$\phi$ by means of a {\em quasi-universal
sheaf}~(\ref{subsec:quSheaf}), as a replacement for the universal
sheaf in the original construction.  (One needs to move to rational
cohomology, instead of the integral one used before.)  While a
universal sheaf may not exist in general, a quasi-universal sheaf
always exists.

\subsection{}
Since Mukai's main interest was to obtain a correspondence between the
cohomology of $M$ and the cohomology of $X$, constructing $\phi$ by
means of a quasi-universal sheaf was enough.  What is lost in this
approach is the equivalence of derived categories -- one can not hope
to get a Fourier-Mukai transform by using a quasi-universal sheaf.  In
this paper we propose a different approach, which makes apparent an
underlying equivalence of derived categories even in the non-fine
case.  Using this idea, we are able to generalize to the case of
non-fine moduli problems the implications $(1)\Rightarrow (2)$ and
$(1)\Rightarrow (3)$ of Orlov's result, and to suggest a conjectural
generalization of the equivalence $(2)\Leftrightarrow (3)$.  The key idea
is that instead of replacing the universal sheaf by a quasi-universal
one, we can replace it by a {\em twisted} universal sheaf.  This
yields an equivalence
\[ \D(M, \alpha) \iso \D(X), \]
where $\alpha$ is the twisting (an element of the Brauer group of $M$,
determined by the original moduli problem data), and $\D(M, \alpha)$
is the derived category of $\alpha$-twisted sheaves on $M$.  (A brief
review of the Brauer group and twisted sheaves is included in
Section~\ref{sec:twsh}.  For further details, the reader is referred
to~\cite{Cal} or~\cite{CalDTw}.)

\subsection{}
We can rephrase this idea in the language of modules over an Azumaya
algebra (which we'll avoid in the sequel, preferring the more
intuitive language of twisted sheaves).  In~\cite{Muk}, Mukai
constructs a module over an Azumaya algebra $\cA$ over $\cO_X$ (i.e.,
a twisted sheaf), but then he forgets the extra structure as an
$\cA$-module, and only uses the $\cO_X$-module structure.  We make use
of the extra structure available to get the equivalence of derived
categories.

\subsection{}
Using twisted sheaves (or modules over an Azumaya algebra) is
particularly relevant in view of current developments of mirror
symmetry.  Recall that one of the fundamental ingredients of
Kontsevich's Homological Mirror Symmetry~(\cite{Kon}) is the derived
category of sheaves on a Calabi-Yau manifold.  Recently there have
been suggestions~(\cite[6.8]{Cal}, \cite{KapOrl}) that in order to get
a mathematical description of the full physical picture, one needs to
study not only derived categories of sheaves, but also derived
categories of twisted sheaves.  In the physical context, the
equivalence $\D(M, \alpha) \iso \D(X)$ suggests that turning on
discrete torsion $\alpha$ on $M$ yields the same physical theory as
having no discrete torsion on $X$.

\subsection{}
Reverting back to questions about moduli spaces of sheaves, it is a
general fact~(\cite[3.3.2 and 3.3.4]{Cal}) that although a universal
sheaf may not exist, there is always a twisting $\alpha\in\Br(M)$ such
that a $\pi_M^*\alpha^{-1}$-twisted universal sheaf exists on $X\times
M$, when $M$ is a moduli space of stable sheaves on $X$.  (Here
$\pi_M$ is the projection from $X\times M$ on the second factor.)  An
$\alpha$ with this property is unique, and it is called the
obstruction to the existence of a universal sheaf on $X\times M$
(because an untwisted universal sheaf exists only if $\alpha=0$).
Once a twisted universal sheaf is found, the proof that it induces an
equivalence of derived categories is almost the same as in the
untwisted case.  What this paper actually brings new is the
calculation of the obstruction $\alpha$ in terms of the map $\phi$
introduced earlier by Mukai.  (This map has a natural interpretation
in the context of twisted sheaves, see Sections~\ref{sec:twsh}
and~\ref{sec:setup}.)

\subsection{}
\label{subsec:prevres}
Before we can state our results, we need to list a few facts about the
situation we study -- $X$ is a K3 surface, $v\in\tH(X, \Z)$ is an
isotropic Mukai vector, and $M$ is the moduli space of stable sheaves
on $X$ whose Mukai vector is $v$ (see~\ref{subsec:MukVect}).  The
vector $v$ and the polarization of $X$ are such that $M$ is non-empty
and compact, so it is a K3 surface (Theorem~\ref{thm:Muk1}).  We do
not assume that $M$ is fine, so Mukai's map we consider is defined on
the rational cohomology, $\phi:\tH(X,\Q)\ra\tH(M,\Q)$.
\begin{enumerate} 
\item The Brauer group $\Br(M)$ of any smooth surface $M$ is
isomorphic to the cohomological Brauer group $\Br'(M)$, defined as
\[ \Br'(M) = H^2_\et(M, \cO_M^*). \]
For a K3 surface there is a natural identification~(\ref{subsec:brqmodz})
\[ \Br'(M) \iso T_M^\chk \otimes \Q/\Z \iso \Hom_\Z(T_M, \Q/\Z). \]
\item The map $\phi$ restricts to an injection $T_X \hookrightarrow 
T_M$, which fits into an exact sequence
\[ 0\ra T_X\stackrel{\phi}{\lra} T_M \ra \Z/n\Z \ra 0, \]
where $n$ is an integer determined by the initial moduli
problem~(\cite[6.4]{Muk}, Theorem~\ref{thm:Muk3}).  Although $\phi$
depends on the choice of quasi-universal sheaf, its restriction to
$T_X$ is independent of this choice.
\item Applying $\Hom_\Z(\scdot, \Q/\Z)$ to the above exact sequence
yields
\[ 0 \ra \Z/n\Z \ra \Hom_\Z(T_M, \Q/\Z) \stackrel{\phi^\chk}{\lra} 
\Hom_\Z(T_X, \Q/\Z) \ra 0, \] 
or, in view of (1), 
\[ 0\ra \Z/n\Z \ra \Br(M) \stackrel{\phi^\chk}{\lra} \Br(X) \ra 0. \]
Elements of $\Ker \phi^\chk$ are those $\alpha\in\Hom_\Z(T_M, \Q/\Z)$
that satisfy $\phi(T_X) \subseteq \Ker \alpha$.
\end{enumerate}

\subsection{}
\label{subsec:defproj}
As a matter of notation, for $w\in \tH(M, \Q)$ the functional 
\[ (w. \scdot)|_{T_M}\bmod\Z \in \Hom_\Z(T_M, \Q/\Z) = \Br(M) \]
will be denoted by $[w]$.  Note that since it is restricted to $T_M$,
$[w]$ only depends on the $H^2(M, \Q)$-component of $w$.  The
condition in (3) above can be written as
\[ [w] \in \Ker\phi^\chk \mbox{ for } w\in \tH(M, \Q) \Leftrightarrow
(w.\phi(t)) \in \Z \mbox{ for all } t\in T_X. \]

\subsection{}
Now we can state our main results:

\begin{theorem}
\label{thm:mainthm1}
Let $X$ be a polarized K3 surface, let $v$ be a primitive isotropic
Mukai vector, and let $M$ be the moduli space of stable sheaves whose
Mukai vector is $v$.  Assume that $M$ is compact and non-empty, and
let $\phi:T_X \ra T_M$ be the restriction of Mukai's map (defined by
means of a quasi-universal or twisted universal sheaf).

If $u\in\tH(X, \Z)$ is such that $(u.v) = 1$ then
$[\phi(u)]\in\Br(M)$ is the obstruction to the existence of a
universal sheaf on $X\times M$.
\end{theorem}

\vbox{
\begin{theorem}
\label{thm:mainthm2}
Under the assumptions of the previous theorem, let $\phi^\chk$ be the
dual of $\phi$, tensored with $\Q/\Z$.  Then we have:
\begin{enumerate}
\item 
the kernel of $\phi^\chk:\Br(M)\ra\Br(X)$ is a cyclic group of order
$n$, generated by the obstruction $\alpha$ to the existence of a
universal sheaf on $X\times M$;
\item 
the map $\phi$ restricts to a Hodge isometry $T_X \iso \Ker \alpha$;
\item 
any $\pi_M^*\alpha^{-1}$-twisted universal sheaf on $X\times M$ 
induces an equivalence of derived categories
\[ \D(M, \alpha) \iso \D(X). \]
\end{enumerate}
\end{theorem}
}

\subsection{}
Theorem~\ref{thm:mainthm2} provides the desired generalization of the
implications $(1)\Rightarrow (2)$ and $(1)\Rightarrow (3)$ in Orlov's
theorem.  For the equivalence $(2)\Leftrightarrow (3)$ we can only
conjecture what the result should be:

\begin{conjecture}
\label{conj:mainconj}
Let $X$ and $Y$ be K3 surfaces, and let $\alpha$ and $\beta$ be
elements in $\Br(X)$ and $\Br(Y)$, respectively.
Using~(\ref{subsec:prevres}.1), $\alpha$ can be identified with a
group homomorphism $T_X\ra \Q/\Z$, and similarly for $\beta$.  Then
the following are equivalent:
\begin{enumerate}
\item the derived categories $\D(X, \alpha)$ and $\D(Y, \beta)$ are
equivalent; 
\item the lattices $\Ker\alpha\subseteq T_X$ and
$\Ker\beta\subseteq T_Y$ are Hodge isometric, where $\Ker\alpha$
and $\Ker\beta$ inherit the Hodge structure from the overlying
lattices $T_X$ and $T_Y$, respectively.
\end{enumerate} 
\end{conjecture}

\subsection{}
This generalizes the original Orlov result, for $\alpha=\beta=0$.  If
$Y$ is a moduli space of stable sheaves on $X$, $\beta\in\Br(Y)$ is
the obstruction to the existence of a universal sheaf on $X\times Y$,
and $\alpha = 0$, the conjecture is a restatement of
Theorem~\ref{thm:mainthm2}: $\phi$ induces a Hodge isometry between
$T_X=\Ker\alpha$ and $\Ker\beta$.

\subsection{}
Conjecture~\ref{conj:mainconj} also generalizes the following result
of Donagi-Pantev, which is relevant in the study of moduli spaces of
twisted Higgs bundles:
\vspace{2mm}

\noindent
{\bf Theorem~(Donagi-Pantev~\cite{DonPan}).}  
{\em 
Let $J/S$ be an elliptic K3 surface with a section, and let $\alpha,
\beta\in \Br(J)$.  Identifying $\Br(J)$ with $\Sha_S(J)$, the
Ogg-Shafarevich group of $J$, yields elliptic K3 surfaces $J_\alpha$
and $J_\beta$ (in general without a section) which correspond to
$\alpha$ and $\beta$, respectively.  Viewing $J$ as a moduli space of
stable sheaves on $J_\alpha$ and $J_\beta$ gives surjections
$\Br(J)\ra \Br(J_\alpha)$ and $\Br(J)\ra\Br(J_\beta)$; let $\bbeta$ be
the image of $\beta$ in $\Br(J_\alpha)$ and $\balpha$ the image of
$\alpha$ in $\Br(J_\beta)$.  Then there exists an equivalence of derived categories of twisted
sheaves $\D(J_\alpha, \bbeta)\iso\D(J_\beta, \balpha)$.}
\vspace{1mm}

A few words of explanation are in order here.  Ogg-Shafarevich theory
associates to an elliptic fibration $X/S$ without a section an
element $\alpha_X$ of the Tate-Shafarevich group $\Sha_S(J)$ of the
relative Jacobian $J/S$ of $X/S$, and this association gives rise to a
bijection between the set of all elliptic fibrations whose relative
Jacobian is $J/S$ and the group $\Sha_S(J)$.  When $J$ is a surface,
$\Sha_S(J)$ is naturally isomorphic to $\Br(J)$.

Given an elliptic fibration $X/S$ without a section, one obtains $J/S$
as the relative moduli space of rank 1, degree 0 semistable sheaves on
the fibers of $X/S$, and thus one gets an obstruction to the existence
of a universal sheaf on $X\times J$ which is an element $\alpha$ in
$\Br(J)$.  It can be shown~(\cite{CalEll}) that $\alpha$ coincides
with the element $\alpha_X\in\Sha_S(J)=\Br(J)$ that classifies $X/S$.

It is not hard to see that Conjecture~\ref{conj:mainconj} implies the
Donagi-Pantev result: we have $T_{J_\alpha} = \Ker \alpha$ and,
viewing elements of the Brauer group as group homomorphisms from the
transcendental lattice to $\Q/\Z$, $\bbeta$ is just the restriction of
$\beta$ to $\Ker \alpha$.  Therefore 
\[ \Ker_{T_{J_\alpha}}\bbeta = \Ker_{T_J} \alpha \cap \Ker_{T_J}
\beta, \]
and similarly,
\[ \Ker_{T_{J_\beta}}\balpha = \Ker_{T_J} \alpha \cap \Ker_{T_J}
\beta, \]
both  these equalities respecting the Hodge structures (being induced on
all the lattices from the Hodge structure of $T_J$).  By
Conjecture~\ref{conj:mainconj},
\[ \D(J_\alpha, \bbeta) \iso \D(J_\beta, \balpha), \] 
which is the Donagi-Pantev result.

In full honesty, the actual statement proven in~\cite{DonPan} is much
stronger: one of $\alpha$ or $\beta$ can be {\em non-algebraic}, i.e.\
an element of $H^2_\an(J, \cO_J^*)$ which is not necessarily torsion.
In this situation the corresponding $J_\alpha$ or $J_\beta$ is
non-algebraic, a situation which we cannot handle.  Therefore
Conjecture~\ref{conj:mainconj} should be thought of as a
generalization of the result in~\cite{DonPan} to the case of algebraic
K3 surfaces.

\subsection{}
In the spirit of our earlier comment~\ref{subsec:sublatcat},
Conjecture~\ref{conj:mainconj} describes the relationship between
sublattices of the Mukai lattice and twisted derived categories.  We
thus get the following dictionary of correspondences between
sublattices of a Mukai lattice endowed with a Hodge structure and
categories on a K3 $X$ ($\alpha\in\Br(X)$):
\begin{center}
\begin{tabular}{c|c}
Lattice         &       Category               \\
\hline
$H^2(X, \Z)$    & $\gCoh(X)$                            \\
$T_X$           & $\D(X)$                               \\
$\Ker(\alpha)$ & $\D(X, \alpha)$        \\
\hline
\end{tabular}
\end{center}

\subsection{}
The paper is organized as follows: in Section~\ref{sec:twsh} we
present general results about the Brauer group on a K3, twisted
sheaves, and derived categories.  We follow up in the next section
with facts about K3 surfaces, along with Mukai's results on moduli
spaces of stable sheaves on them.  Section~\ref{sec:defos} deals with
deformations of vector bundles as twisted sheaves, and in
Section~\ref{sec:mainthm} we prove the main theorem and discuss some
possible ways of approaching the proof of
Conjecture~\ref{conj:mainconj}.
\vspace{2mm}

\noindent\textbf{Conventions.}
All our spaces are complex manifolds over $\C$, and the topology used
is either the analytic or \'etale one.  When referring to derived
categories, we mean the bounded derived category of complexes with
coherent cohomology.
\vspace{2mm}

\noindent\textbf{Acknowledgments.}  
The results in this paper are part of my Ph.D.\ work, completed at
Cornell University.  I would like to thank my supervisor, Mark Gross,
for taking the time to teach me algebraic geometry, and for providing
plenty of help and encouragement.  It is a pleasure also to thank Tony
Pantev, discussions with whom have been extremely useful in clarifying
some of the main ideas in this work.

\section{Twisted sheaves and derived categories}
\label{sec:twsh}

In this section we list some results about the Brauer group of a
scheme.  Our reference for this topic is~\cite[Chapter IV]{Mil}.  We
also include a sketch of the definition and main properties of twisted
sheaves.  The reader unfamiliar with the subject is referred
to~\cite[Chapters 1 and 2]{Cal} or~\cite{CalDTw}.  The topology used,
unless otherwise mentioned, is the \'etale or analytic topology.

\subsection{}
The {\em cohomological Brauer group} of a scheme $X$ is the group
\[ \Br'(X) = H^2_\et(X, \cO_X^*). \]
It naturally occurs in many aspects of algebraic geometry, as a higher
generalization of the Picard group.  To have an example in mind
consider the question of classifying projective bundles over a space
$X$, up to those bundles that are projectivizations of vector bundles.
In the \'etale topology the sequence of sheaves of groups
\[ 0 \ra \cO_X^* \ra \GL(n) \ra \PGL(n) \ra 0, \]
is exact, hence it yields the exact sequence
\[ H^1(X, \GL(n)) \ra H^1(X, \PGL(n)) \ra H^2(X, \cO_X^*) = \Br'(X). \]
We read this as saying that the obstruction to lifting a projective
bundle (given by an element $p$ of $H^1(X, \PGL(n))$) to a vector bundle
(element of $H^1(X, \GL(n))$) is the image of $p$ in the cohomological 
Brauer group, under the coboundary map.

\subsection{}
In fact, we are more interested in a subgroup of $\Br'(X)$, namely the
image of the coboundary maps $H^1(X, \PGL(n)) \ra \Br'(X)$ for all
$n$.  This subgroup is the {\em Brauer group} of $X$, denoted by
$\Br(X)$.  We list below some of its main properties.

The Brauer group is torsion; this follows from the short exact
sequence
\[ 0 \ra \Z/n\Z \ra \SL(n) \ra \PGL(n) \ra 0, \]
by taking the cohomology long exact sequence, and deducing that the
image of $\PGL(n)$ in $\Br'(X)$ is contained in the image of the map
$H^2(X, \Z/n\Z) \ra H^2(X, \cO_X^*)$, and hence is $n$-torsion.

If $X$ is smooth, $\Br'(X)$ is torsion as well~(\cite[II, 1.4]{Mil}).
If $X$ is a smooth curve, $\Br(X) = \Br'(X) = 0$.  If $X$ is a smooth
surface, $\Br(X) = \Br'(X)$~(\cite[IV, 2.16]{Mil}).

\subsection{}
Consider the Kummer sequence
\[ 0 \ra \Z/n\Z \ra \cO_X^* \stackrel{\scdot^n}{\lra} \cO_X^* \ra 0, \]
which is exact in both the \'etale and analytic topologies.  Taking
the associated long exact sequence yields
\[ 
\begin{diagram}
\Pic(X) & \rTo^{\scriptstyle{\cdot n}} & \Pic(X) & \rTo^{\scriptstyle{c_1 \bmod
n}} & H^2(X, \Z/n\Z) & \rTo & \Br'(X) & \rTo^{\scriptstyle{\cdot n}} & \Br'(X),
\end{diagram}
\]
which implies that the $n$-torsion part of $\Br'(X)$, $\Br'(X)_n$,
fits in the exact sequence
\[ 0 \ra \Pic(X) \otimes \Z/n\Z \ra H^2(X, \Z/n\Z) \ra \Br'(X)_n \ra 0. \]
Taking the direct limit over all $n$, we conclude that on any scheme or
analytic space $X$ we have the exact sequence
\[ 0 \ra \Pic(X)\otimes \Q/\Z \ra H^2(X, \Q/\Z) \ra \Br'(X)_\tors \ra 0. \]
If $X$ is a smooth scheme over $\C$, and $X^\an$ is the associated
analytic space, then we have
\[ \Br'(X) = \Br'_\an(X)_\tors, \]
because $\Pic(X)$ and $H^2(X, \Q/\Z)$ are the same in the \'etale and
analytic topologies.

\subsection{}
\label{subsec:brqmodz}
Specializing to the case of a K3 surface $X$, we have $H^2(X, \Q/\Z)
\iso H^2(X, \Z) \otimes \Q/\Z$ (because $H^3(X, \Z) = 0$), and
hence
\[ \Br(X) \iso (H^2(X, \Z) / \NS(X)) \otimes \Q/\Z. \]
There is a natural isomorphism 
\[ H^2(X, \Z) / \NS(X) \iso T_X^\chk, \]
which maps $v\in H^2(X, \Z) / \NS(X)$ to the functional $(v,
\scdot)$ restricted to $T_X$.  (To prove that this map is an
isomorphism one needs to use the fact that $T_X$ is a primitive
sublattice of the unimodular lattice $H^2(X, \Z)$.)  We conclude that
on any K3 surface $X$ there is a natural isomorphism
\[ \Br(X) \iso T_X^\chk \otimes \Q/\Z = \Hom_\Z(T_X, \Q/\Z). \]

\subsection{}
We now shift our attention to the topic of twisted sheaves.  Let $X$
be a scheme or analytic space, and let $\alpha\in \Br'(X) = H^2(X,
\cO_X^*)$ be represented by a \Cech 2-cocycle, given along a fixed
open cover $\{U_i\}_{i\in I}$ by sections
\[ \alpha_{ijk}\in \Gamma(U_i\cap U_j\cap U_k, \cO_X^*). \]
An $\alpha$-twisted sheaf $\cF$ (along the fixed cover) consists of a
pair 
\[ (\{\cF_i\}_{i\in I}, \{\phi_{ij}\}_{i,j\in I}), \]
where $\cF_i$ is a sheaf on $U_i$ for all $i\in I$ and
\[ \phi_{ij}: \cF_j|_{U_i\cap U_j} \ra \cF_i|_{U_i\cap U_j} \]
is an isomorphism for all $i,j\in I$, subject to the conditions:
\begin{enumerate}
\item $\phi_{ii} = \id$;
\item $\phi_{ij} = \phi_{ji}^{-1}$;
\item $\phi_{ij}\circ \phi_{jk}\circ \phi_{ki} = \alpha_{ijk} \cdot \id$.
\end{enumerate}

The class of $\alpha$-twisted sheaves together with the obvious notion of
homomorphism is an abelian category, denoted by $\gMod(X, \alpha)$,
the category of $\alpha$-twisted sheaves.  If one requires all the
sheaves $\cF_i$ to be coherent, one obtains the category of coherent
$\alpha$-twisted sheaves, denoted by $\gCoh(X, \alpha)$.

This notation is consistent, since one can prove that these categories
are independent of the choice of the covering
$\{U_i\}$~(\cite[1.2.3]{Cal}) or of the particular cocycle
$\{\alpha_{ijk}\}$~(\cite[1.2.8]{Cal}) (all the resulting categories
are equivalent to one another).

\subsection{}
\label{subsec:functors}
For $\cF$ an $\alpha$-twisted sheaf, and $\cG$ an $\alpha'$-twisted
sheaf, one can define $\cF\otimes \cG$ (which is an
$\alpha\alpha'$-twisted sheaf), as well as $\sHom(\cF, \cG)$ (which is
$\alpha^{-1}\alpha'$-twisted), by gluing together the corresponding
sheaves.  If $f:Y\ra X$ is any morphism, $f^*\cF$ is an
$f^*\alpha$-twisted sheaf on $Y$.  Finally, if $\cF\in \gMod(Y,
f^*\alpha)$, one can define $f_* \cF$, which is $\alpha$-twisted on
$X$.  It is important to note here that one can not define arbitrary
push-forwards of twisted sheaves.

These operations satisfy all the usual relations (adjointness of $f_*$
and $f^*$, relations between $\sHom$ and $\otimes$, etc.)

The category $\gMod(X, \alpha)$ has enough injectives, and enough
$\cO_X$-flats~(\cite[2.1.1, 2.1.2]{Cal}).

\subsection{}
We are mainly interested in $\D(\gMod(X, \alpha))$, the derived
category of complexes of $\alpha$-twisted sheaves on $X$ with coherent
cohomology.  For brevity, we'll denoted it by $\D(X, \alpha)$.  Since
the category $\gCoh(X, \alpha)$ does not have locally free sheaves of
finite rank if $\alpha\not\in\Br(X)$, from here on we'll only consider
the case $\alpha\in\Br(X)$.

The technical details of the inner workings of $\D(X, \alpha)$ can be
found in~\cite{CalDTw} or~\cite[Chapter 2]{Cal}.  The important facts
are that one can define derived functors for all the functors
considered in~\ref{subsec:functors}, and they satisfy the same
relations as the untwisted ones (see for example~\cite[II.5]{HarRD}).
One can prove duality for a smooth morphism $f:X\ra Y$, which provides
a right adjoint
\[ f^!(\scdot) = \Ld f^*(\scdot) \otimes_X \omega_{X/Y}[\dim_X Y] \]
to $\R f_*(\scdot)$, as functors between $\D(Y, \alpha)$ and $\D(X,
f^* \alpha)$.  

\subsection{}
If $X$ and $Y$ are smooth schemes or analytic spaces,
$\alpha\in\Br(Y)$, and $\dcE\in\D(X\times Y, \pi_Y^*\alpha^{-1})$
(where $\pi_X$ and $\pi_Y$ are the projections from $X\times Y$ to $X$
and $Y$ respectively), we define the integral functor
\begin{align*}
\FMYX^\dcE & :\D(Y, \alpha) \ra \D(X), 
\intertext{given by}
\FMYX^\dcE(\scdot) & = \pi_{X,*}(\pi_Y^*(\scdot) \lotimes \dcE).
\end{align*}

The following criterion for determining when $\FMYX^{\dcE}$ is an
equivalence (whose proof can be found in~\cite{CalDTw}
or~\cite[3.2.1]{Cal}) is entirely similar to the corresponding ones
for untwisted derived categories due to Mukai~\cite{MukAb},
Bondal-Orlov~\cite{BonOrl} and Bridgeland~\cite{Bri}.

\begin{theorem}
\label{thm:equiv}
The functor $F=\FMYX^\dcE$ is fully faithful if and only if for each
point $y\in Y$,
\[ \Hom_{\D(X)}(F\cO_y, F\cO_y) = \C, \]
and for each pair of points $y_1, y_2\in Y$, and each integer $i$,
\[ \Ext^i_{\D(X)}(F\cO_{y_1}, F\cO_{y_2}) = 0 \]
unless $y_1=y_2$ and $0\leq i\leq \dim Y$.  (Here $\cO_y$ is the
skyscraper sheaf $\C$ on $y$, which is naturally an $\alpha$-sheaf.)

Assuming the above conditions satisfied, then $F$ is an equivalence of
categories if and only if for every point $y\in Y$,
\[ F\cO_y \lotimes \omega_X \iso F\cO_y.\]
\end{theorem}

\subsection{}
\label{subsec:chclass}
Since we want to relate derived categories to cohomology, we'd like to
define the notion of Chern character for twisted sheaves on a space
$X$.  We do this as follows: we fix, once and for all, a locally free
$\alpha^{-1}$-twisted sheaf $\cE$, and define the Chern character
of an $\alpha$-twisted sheaf $\cF$ to be
\[ \ch_\cE(\cF) = \frac{1}{\rk(\cE)} \ch(\cF\otimes \cE), \]
where the right hand side of the equality is computed as the Chern
character of the usual (untwisted) sheaf $\cF\otimes\cE$.  Note that
if we define Chern classes in this way, they will live in the rational
cohomology of $X$, not in the integral cohomology as the usual ones.
The factor $1/\rk(\cE)$ is introduced so that the Chern character of a
point is the expected one.

This definition is dependent on the choice of $\cE$, and therefore it
is important to find out how our Chern character changes when using
different $\cE$'s.  Given any two locally free $\alpha$-twisted
sheaves $\cE$ and $\cE'$, there exist locally free, untwisted sheaves
$\cG$, $\cG'$ such that $\cE\otimes \cG\iso \cE'\otimes\cG'$ (for
example, one can take $\cG = \cE^\chk \otimes \cE'$ and $\cG' =
\cE^\chk\otimes \cE$), so it is enough to assume that $\cE' = \cE\otimes
\cG'$ for a locally free, untwisted $\cG'$.  Then an easy calculation
shows that
\[ \ch_{\cE'}(\cF) = \ch_\cE(\cF).\frac{\ch(\cG')}{\rk(\cG')} \]
(see also~\cite[Proof of Theorem 1.4 and 1.5, p. 385]{Muk}).  We'll
see~(\ref{subsec:phitrindep}) that this implies that when we define
maps on transcendental parts of the cohomology of $X$ using the above
definition of the Chern character for twisted sheaves, the choice of
$\cE$ is irrelevant.

\section{Mukai's results}
\label{sec:setup}

In order to fix the notation, we present in this section a few results
about K3 surfaces and moduli spaces of stable sheaves on them.  Most
of these results are either classical, or are due to Mukai~\cite{Muk},
although we present them using the language of twisted sheaves.

\subsection{}
Mukai's main idea in~\cite{Muk} is to use a universal sheaf (possibly
twisted, or a quasi-universal one) to define a map between sheaves on
$X$ and sheaves on $M$ (more precisely, between the derived categories
of $X$ and $M$).  Taking a modified version of the Chern character
yields a map between the algebraic parts of the cohomology of $X$ and
of $M$, which can be extended to the full cohomology.  Since the map
on derived categories was an equivalence, the map on cohomology is an
isomorphism, even after being extended to the total cohomology groups
of $X$ and $M$.  Furthermore, since an equivalence of categories
preserves the relative Euler characteristic of two complexes of
sheaves (the alternating sum of the dimensions of the Ext groups),
this gives a bilinear form that is preserved by the map on cohomology.
We expand these ideas in the following few paragraphs.

\subsection{}
Let $X$ be a complex K3 surface, in other words, a compact complex
manifold of complex dimension 2, simply connected and with $K_X=0$.
We have $H^2(X, \Z) = \Z^{22}$, and considering this group with the
intersection pairing we obtain a lattice which is isomorphic to
\[ \sL_\Kt = E_8^{\oplus 2} \oplus U(1)^{\oplus 3}. \]
Inside the $H^2(X, \Z)$ lattice there are two natural sublattices, the
N\'eron-Severi sublattice of $X$, $\NS(X)$, (consisting of first Chern
classes of holomorphic vector bundles), and its orthogonal complement,
the transcendental lattice $T_X = \NS(X)^\perp$.  Both these lattices
are primitive sublattices of $H^2(X, \Z)$, but may be non-unimodular.

The complex structure of $X$ is reflected in the Hodge decomposition
of $H^2(X, \C)$,
\[ H^2(X, \C) = H^{2,0}(X) \oplus H^{1,1}(X) \oplus H^{0,2}(X), \]
and these groups are 1-, 20-, and 1-dimensional, respectively.  This
decomposition induces in turn a Hodge structure on $T_X$ (since
$H^{2,0}(X)$ is orthogonal to any algebraic vector).

Two lattices $\sL, \sL'$, endowed with Hodge structures, will be said
to be {\em Hodge isometric} if there is an isometry between them,
preserving the Hodge structure.  As an application of this concept,
the Torelli theorem can be stated by saying that two K3 surfaces $X$
and $Y$ are isomorphic if and only if $H^2(X, \Z)$ and $H^2(Y, \Z)$
are Hodge isometric.

\subsection{}
\label{subsec:MukProd}
The {\em Mukai lattice} of $X$ is defined to be
\[ \tH(X, \Z) = H^0(X, \Z) \oplus H^2(X, \Z) \oplus H^4(X, \Z), \]
endowed with the product
\[ ((r,l,s).(r',l',s')) = \int_X l.l' - r.s' - r'.s, \]
where the dot product on the right hand side is the cup product in
$H^*(X, \Z)$.

\subsection{}
\label{subsec:MukHodge}
The {\em Hodge decomposition} on $\tH(X, \Z)$ is given by 
\begin{align*}
\tH^{2,0}(X) & = H^{2,0}(X) \\
\tH^{1,1}(X) & = H^0(X, \C) \oplus H^{1,1}(X) \oplus H^4(X, \C) \\
\tH^{0,2}(X) & = H^{0,2}(X)
\end{align*}
Elements in $\tH^{1,1}(X)$ will be called {\em algebraic}.

We'll sometimes also consider $\tH(X, \Q) = \tH(X, \Z) \otimes \Q$, with
intersection pro\-duct and Hodge decomposition defined in a similar fashion.

\subsection{}
\label{subsec:MukVect}
For a coherent sheaf $\cF$ (or, more generally, an element of
$\D(X)$), define $v(\cF) \in \tH(X, \Z)$ by
\[ v(\cF) = \ch(\cF).\sqrt{\Td(X)} = (\rk(\cF), c_1(\cF), \rk(\cF)\omega +
\frac{1}{2}c_1(\cF)^2 - c_2(\cF)), \] 
(where $\omega\in H^4(X, \Z)$ is the fundamental class of $X$).  This
element is called the {\em Mukai vector} of $\cE$.  From
Grothendieck-Riemann-Roch it follows that
\[ \chi(\cE, \cF) = \sum (-1)^i \dim\Ext_X^i(\cE, \cF) = -(v(\cE).
v(\cF)), \]
for any $\cE, \cF\in\D(X)$.

If $\cF$ is an $\alpha$-twisted sheaf on $M$, we define its Mukai
vector using the same formula as in the untwisted case, but using the
definition of the Chern character in~\ref{subsec:chclass}.  Since this
definition depends upon the choice of a locally free,
$\alpha^{-1}$-twisted sheaf $\cE$ on $M$, we'll denote this type of
Mukai vector by $v_\cE(\cF)$.

\begin{theorem}[{\cite[Theorem 1.4]{Muk}}]
\label{thm:Muk1}
Let $X$ be a polarized K3 surface, and let $v\in\tH(X, \Z)$ be a
primitive (indivisible) vector that lies in the algebraic part of
$\tH(X)$.  Assume that $v$ is isotropic (i.e.\ $(v.v) = 0$), and that
the moduli space of stable sheaves of Mukai vector $v$, $M(v)$, is
non-empty and compact.  Then $M(v)$ is a K3 surface.
\end{theorem}

\subsection{}
\label{subsec:quSheaf}
Under the assumptions of the previous theorem, let $M=M(v)$.  There
exists a unique element $\alpha\in\Br(M)$, such that a
$\pi_M^*\alpha^{-1}$-twisted universal sheaf $\cP$ exists on $X\times
M$~(\cite[3.3.2 and 3.3.4]{Cal}).  The twisting $\alpha$ is called the
{\em obstruction to the existence of a universal sheaf} on $X\times
M$.  Let $\cE$ be a fixed $\alpha$-twisted locally free sheaf of
finite rank on $M$, (which exists by~\cite[1.3.5]{Cal}), and let
$\tilde{\cP\ }= \cP \otimes \pi_M^* \cE$.  It is an untwisted sheaf on
$X\times M$, and any such sheaf will be called a {\em quasi-universal
sheaf}.  It has the property that
\[ \tilde{\cP|}_{X\times [\cF]} \iso \cF^{\oplus n} \]
for some $n$ that only depends on $\tilde{\cP.}$\, (Here, $\cF$ is a
stable sheaf on $X$ with Mukai vector $v$, and $[\cF]$ is the point of
$M$ that corresponds to it.)  In fact, $\tilde{\cP\ }$satisfies a
certain universal property which is very similar to that enjoyed by a
universal sheaf; see~\cite[Appendix 2]{Muk}.  Mukai uses a
quasi-universal sheaf to define the correspondence between the
cohomology of $X$ and $M$, however, since it seems more natural, we'll
avoid using this quasi-universal sheaf and use the twisted universal
sheaf $\cP$ instead.

\subsection{}
\label{subsec:MukIsom}
The dual $\cQ$ of $\cP$ is defined by
\[ \cQ = \R\sHom(\cP, \cO_{X\times M}), \]
as an element of $\D(X\times M, \pi_M^*\alpha)$.  Using it, we can
define the correspondence
\begin{align*}
\phi = \fmXM^{\cQ, \cE} & : \tH(X, \Q) \ra \tH(M, \Q) 
\intertext{given by the formula}
\phi_\cE(\scdot) & = \pi_{M,*}(\pi_X^*(\scdot).v_{\pi_M^*\cE}(\cQ)),
\intertext{where}
v_{\pi_M^*\cE}(\cQ) & = \ch_{\pi_M^*\cE}(\cQ).\sqrt{\Td(X\times M)},
\end{align*}
and $\ch_{\pi_M^*\cE}(\cQ)$ is the Chern character of the twisted
sheaf $\cQ$, defined in~\ref{subsec:chclass}, computed using the
$\alpha^{-1}$-twisted locally free sheaf $\cE$ on $M$.

The reason one uses $\cQ$ is the fact that if $\FMMX^\cP$ is the
integral transform associated to $\cP$ (which is an equivalence), then
$\FMXM^\cQ$ is its inverse.  Also, from Grothendieck-Riemann-Roch it
follows that
\[ v_\cE(\FMXM^\cQ(\cF)) = \fmXM^{\cQ, \cE}(v(\cF)). \]

\subsection{}
This is the core of the following result:
\begin{theorem}[{\cite[Theorem 1.5]{Muk}}]
\label{thm:isomK3}
Under the hypotheses of Theorem~\ref{thm:Muk1}, the map $\phi$ is a
Hodge isometry between $\tH(X, \Q)$ and $\tH(M, \Q)$.  It maps $v\in
\tH(X, \Q)$ to the vector $(0,0,\omega)\in \tH(M, \Q)$, and it
therefore induces a Hodge isometry
\[ v^\perp/v \iso H^2(M, \Q), \]
the former computed inside $\tH(X, \Q)$.  Restricted to $v^\perp$,
this isometry is independent of the choice of $\cE$, and is integral,
i.e.\ it takes integral vectors to integral vectors.  It therefore
induces a Hodge isometry
\[ H^2(M, \Z) \iso v^\perp/v, \]
the latter now computed in $\tH(X, \Z)$.  
\end{theorem}

\begin{remark}
Using the Torelli theorem, this gives a complete description of the moduli
space $M$.
\end{remark}

\begin{theorem}
\label{thm:Muk3}
Let $n$ be the greatest common divisor of the numbers $(u.v)$, where
$u$ runs over all $\tH^{1,1}(X)\cap \tH(X,\Z)$ ($v$ is the Mukai
vector referred to in Theorem~\ref{thm:Muk1}).  Then the following
statements hold:
\begin{enumerate}
\item There exist $\alpha$-twisted locally free sheaves on $M$ of
ranks $r_1, r_2, \ldots, r_k$, with $\gcd(r_1, r_2, \ldots, r_k) = n$,
and therefore $\alpha$ is $n$-torsion.  

\item Any map $\phi = \fmXM^{\cQ, \cE}$ maps $T_X$ into $T_M$ (viewing
$T_X$ as a sublattice of $\tH(X, \Z)$ via the inclusion $\lambda
\mapsto (0,\lambda, 0)$), and the restriction $\phi|_{T_X}$ is
independent of the choice of $\cE$.

\item There exists $\lambda\in T_X$ such that $v-\lambda$ is divisible 
by $n$ (in $\tH(X, \Z)$); for such a $\lambda$ we have $\phi(\lambda)$
divisible by $n$ (in $T_M$).

\item $\phi|_{T_X}$ is injective, and its cokernel is a finite, cyclic
group of order $n$, generated by $\phi(\lambda)/n$ for any $\lambda$
satisfying the condition in (3).
\end{enumerate}
\end{theorem}

\begin{proof}
For the first statement, see~\cite[Remark A.7]{Muk}.  The other
statements are~\cite[6.4]{Muk}.  
\end{proof}

\subsection{}
\label{subsec:phitrindep}
The calculation in~\ref{subsec:chclass} shows that for any $u\in
\tH(X, \Q)$, changing $\cE$ will only change $\phi(u)$ by an algebraic
amount, and therefore $[\phi(u)]$ (as defined in~\ref{subsec:defproj})
is independent of the choice of $\cE$.

\section{Twisted deformations of vector bundles}
\label{sec:defos}

In this section we study what happens when we try to deform a vector
bundle from the central fiber of a family, if the first Chern class of
the vector bundle fails to deform to neighboring fibers.  We show that
if such a deformation exists as a twisted sheaf, then there is a
simple formula for what the twisting needs to be.  This will be used
in the proof of Theorem~\ref{thm:mainthm1} to identify the obstruction
to the existence of a universal sheaf by a deformation argument, but
the results in this section may be of independent interest.

\subsection{}
\label{subsec:defass}
Let's first set up the context.  We start with $f:X\ra S$, a proper,
smooth morphism of analytic spaces, and with $0$ a closed point of
$S$.  The Brauer group we consider is $\Br'_\an(X)_\tors$, which is
the natural generalization to the analytic setting of the \'etale
Brauer group used in the algebraic case.  Throughout this section
we'll be loose in our notation and refer to $\Br'_\an(X)_\tors$ as the
Brauer group of $X$, or $\Br'(X)$.

Let $X_0$ be the fiber of $f$ over $0$.  We consider an element
$\alpha\in\Br'(X)$, such that $\alpha|_{X_0}$ is trivial as an element
of $\Br'(X_0)$, and we assume we are given a locally free
$\alpha$-twisted sheaf $\cE$ on $X$.  Since $\alpha|_{X_0} = 0$, we
can modify the transition functions of $\cE|_{X_0}$ by a coboundary in
such a way that we get an untwisted locally free sheaf $\cE_0$ on
$X_0$.  We want to understand what happens to $c_1(\cE_0)$ in the
neighboring fibers.  The actual value of $c_1(\cE_0)$ depends on the
choice of coboundary used to trivialize $\alpha|_{X_0}$ ($\cE_0$ could
change by the twist by a line bundle), but the image of
$c_1(\cE_0)/\rk(\cE_0)$ in $H^2(X_0, \Z)\otimes \Q/\Z$ is independent of
these choices.

Since the morphism $f$ is smooth, by possibly restricting first the
base $S$ to a smaller, simply connected one, the restriction
homomorphisms provide identifications
\[ H^i(X, \Z) \iso H^i(X_t, \Z) \]
for all $i\geq 0$ and all $t\in S$.

For any space $X$ we have
\[ H^2(X, \Q/\Z) = H^2(X, \Z) \otimes \Q/\Z \oplus H^3(X, \Z)_\tors, \]
from the universal coefficient theorem.  We saw~(\ref{subsec:brqmodz}) that
$\Br'(X)$ is the quotient of $H^2(X, \Q/\Z)$ by the image of $\Pic(X)\otimes
\Q/\Z$.  But classes in $H^3(X, \Z)_\tors$ cannot become zero in this quotient,
so (in view of the fact that cohomology groups of the fibers are locally
constant over the base) the only way an element of $\Br'(X)$ can become zero in
a central fiber $X_0$ without being zero in the neighboring fibers is if it
belongs to $H^2(X, \Z)\otimes \Q/\Z$, and in the central fiber it is also in
the image of $\Pic(X_0)\otimes \Q/\Z$.

\subsection{}
For an element $w\in H^2(X, \Q)$, we'll denote by $[w]$ the image of
$w$ in 
\[ (H^2(X, \Z)/\NS(X)) \otimes \Q/\Z \subseteq \Br'(X). \]
Because of the considerations in~\ref{subsec:defass}, we can write
$\alpha$ as $[c/n]$, for some class $c\in H^2(X, \Z)$ and $n\in \Z$.
The fact that $\alpha|_{X_0} = 0$ means that $c|_{X_0}$ belongs to
$\Pic(X_0)$.  Our goal is to identify $c|_{X_0}$ and $n$, in terms of
the locally free sheaf $\cE_0$.  This is the content of the following
theorem:

\begin{theorem}
\label{thm:localpha}
Let $\cE$ be a locally free $\alpha$-twisted sheaf on $X$, and let
$\cE_0 = \cE|_{X_0}$.  Assume that $S$ is small enough (say,
contractible), so that we have an identification $H^i(X, \Z) \iso
H^i(X_t, \Z)$ for all $i$ and all $t\in S$.  We assume that
$\alpha|_{X_0} = 0$, and therefore we can modify the transition
functions of $\cE_0$ so that it is an usual sheaf on $X_0$.  Then we
have
\[ \alpha = \left [ -\frac{1}{\rk(\cE_0)}c_1(\cE_0) \right ]. \]
\end{theorem}

The interpretation of this theorem is the following: if we try to
deform a vector bundle $\cE_0$, given on the central fiber $X_0$, in a
family in which the class $c_1(\cE_0)$ is not algebraic in neighboring
fibers $X_t$, the only hope to be able to do this is to deform $\cE_0$
as a twisted sheaf, and then the twisting should be precisely
\[ \left [ -\frac{1}{\rk(\cE_0)}c_1(\cE_0) \right ]. \]

\subsection{}
The idea of the proof is quite straightforward: given a locally free
sheaf $\cE$ (twisted or not) on a space $X$, we consider its associated
projective bundle.  For any projective bundle we define a topological
invariant (the {\em topological twisting characteristic}), which is an
element of $H^2(X, \Z/n\Z)$ (where $n = \rk(\cE)$).  This topological
invariant behaves well with respect to restriction, and it is related
to $c_1(\cE)/\rk(\cE)$ when $\cE$ is not twisted, and to $\alpha$ when
$\cE$ is $\alpha$-twisted.  This enables us to compare $\alpha$ and
$c_1(\cE_0)/\rk(\cE_0)$ in the situation we are interested in.

\subsection{}
Consider the two short exact sequences:
\begin{align*}
0 \ra \cO_X^* \ra & \GL(n) \ra \PGL(n) \ra 0
\intertext{and}
0 \ra \Z/n\Z \ra & \SL(n) \ra \PGL(n) \ra 0.
\end{align*}
For an element $p$ of $H^1(X, \PGL(n))$, let $a(p)$ and $t(p)$ be the
images of $p$ under the two coboundary maps
\begin{align*}
H^1(X, \PGL(n)) & \ra H^2(X, \cO_X^*)
\intertext{and}
H^1(X, \PGL(n)) & \ra H^2(X, \Z/n\Z),
\end{align*}
respectively.  We'll call $a(p)$ the {\em analytic twisting class} of
$p$ and $t(p)$ the {\em topological twisting class} of $p$.  The first
one belongs to $\Br'(X)$, and as such depends on the complex structure
of $X$, while the second one is in $H^2(X, \Z/n\Z)$ and depends only
on the topology of $X$.

If $Y\ra X$ is a $\pj^{n-1}$-bundle over $X$, then the analytic and
topological twisting classes of $Y/X$, $t(Y/X)$ and $a(Y/X)$, are the
classes of the element of $H^1(X, \PGL(n))$ associated to the bundle
$Y\ra X$.

Note that if $\cE$ is an $\alpha$-twisted locally free sheaf
on $X$, we can consider its associated projective bundle $Y\ra X$
(which makes sense even in the twisted case), and then its analytic
twisting class satisfies
\[ a(Y/X) = \alpha. \]

\subsection{}
The following proposition computes $t(Y/X)$ when $Y/X$ is the
projectivization of a locally free (untwisted) sheaf of rank $n$ on
$X$:

\begin{proposition}
\label{t1En}
Let $X$ be a scheme or analytic space, $\cE$ a rank $n$ locally free
sheaf on $X$, and let $Y=\bProj(\cE)\ra X$ be the associated projective
bundle.  Then
\[ t(Y/X) = -c_1(\cE) \bmod n. \]
(Here, and in the sequel, by reducing mod $n$ we mean applying the natural map
$H^2(X, \Z) \ra H^2(X, \Z/n\Z)$).
\end{proposition}

\begin{proof}
Consider the commutative diagram with exact rows and diagonals
\[
\begin{diagram}[silent,height=1.5em,width=2.25em,tight,dpi=600]
	& 	&   0	&	&	&	&   		&		&	&	&   0	\\
	&	&	& \rdTo	&	& 	&		&		&	& \ruTo	&	\\
   0	& \rTo	& \Z/n\Z	& \rTo	&\SL(n)	& 	& \rTo		&		&\PGL(n)& \rTo	&   0	\\
	&	&	&	&	& \rdTo	&		& \ruTo		&	&	&	\\
	&	&\dEqual&	&	&	&\GL(n)		&		&	&	&	\\
	&	&	&	&	& \ruTo &		& \rdTo^{\det}	&	&	&	\\
   0	& \rTo	& \Z/n\Z	& \rTo	& \cO_X^*	&	&\rTo^{\cdot^n}	&		& \cO_X^*	& \rTo	&   0	\\
	&	&	& \ruTo	&	&	&		&		&	& \rdTo	&   	\\
	& 	&   0	&	&	&	&   		&		&	&	&   0.
\end{diagram}
\]
By an easy exercise in homological algebra we get the anti-commutative
diagram
\[
\begin{diagram}[height=1.5em,width=2em,silent,dpi=600]
		&		& H^1(X, \PGL(n))	& \rTo^{t_1}		& H^2(X, \Z/n\Z)	\\
		& \ruTo 	&			&			&		\\
H^1(X, \GL(n))	&		&			&			& \dEqual	\\
		& \rdTo^{\det} 	&			&			&		\\
		&		& H^1(X, \cO_X^*)		& \rTo^{c_1 \bmod n}	& H^2(X, \Z/n\Z)
\end{diagram}
\]
which is precisely what we need.
\end{proof} 

\begin{proposition}
\label{compantop}
For any integer $n$, the following diagram commutes
\[
\begin{diagram}[height=2em,width=2em,dpi=600,silent]
H^2(X, \Z)	& 	& \rTo^{\mod n}	& 	& H^2(X, \Z/n\Z) \\
		& \rdTo_{x\mapsto\left[\frac{1}{n}x\right]} &	& \ldTo_{p} & \\
		&	& \Br'(X),&	&		
\end{diagram}
\]
where the map $p:H^2(X, \Z/n\Z) \ra \Br'(X)$ is obtained from the
natural inclusion $\Z/n\Z \hookrightarrow \cO_X^*$, and the map
$H^2(X, \Z) \ra \Br'(X)$ is taking $x\in H^2(X, \Z)$ to
$\frac{1}{n}x\in H^2(X, \Z) \otimes \Q/\Z$, which then maps to
$\Br'(X)$.

Furthermore, if $Y\ra X$ is a $\pj^{n-1}$-bundle over $X$, we have
\[ p(t(Y/X)) = a(Y/X). \]
\end{proposition}

\begin{proof}
Trivial chase through the definitions.  The last statement follows from the
commutativity of the diagram
\[ 
\begin{diagram}[height=2em,dpi=600]
H^1(X, PGL(n)) 	& \rTo^{t}	& H^2(X, \Z/n\Z)	\\
\dEqual		&		& \dTo_p		\\
H^1(X, PGL(n))	& \rTo^{a}	& \Br'(X),
\end{diagram}
\]
which is deduced from the map of short exact sequences
\[
\begin{diagram}[height=2em,dpi=600]
0 & \rTo & \Z/n\Z  & \rTo & \SL(n) & \rTo & \PGL(n) & \rTo & 0 \\
  &      & \dTo    &      & \dTo   &      & \dEqual &      &   \\
0 & \rTo & \cO_X^* & \rTo & \GL(n) & \rTo & \PGL(n) & \rTo & 0.
\end{diagram}
\]
\end{proof}

\begin{proof}[Proof of Theorem~\ref{thm:localpha}]
Let $n=\rk(\cE)$, and consider the projective bundle associated to $\cE$,
$Y=\bProj(\cE)\ra X$.  Using the naturality of the topological twisting
class we get
\[ \alpha = a(Y/X) = p(t(Y/X)) = p(t(Y_0/X_0)) =
p(-c_1(\cE_0)\bmod n) = \left[-\frac{1}{n} c_1(\cE_0)\right], \] 
where the last three equalities are to be understood via the
identification 
\[ H^2(X,\Z/n\Z)\iso H^2(X_0, \Z/n\Z), \]
in other words $t(Y_0/X_0)$ and $c_1(\cE_0)$ are considered as classes
in $H^2(X,\Z/n\Z)$ and $H^2(X, \Z)$, respectively.
\end{proof}

\section{The proof of the main theorems}
\label{sec:mainthm}

\subsection{}
We consider again the setup of Theorem~\ref{thm:Muk1}: $X$ is a K3
surface, $v$ is a primitive, isotropic Mukai vector on $X$, and $M$ is
the moduli space of stable sheaves on $X$ whose Mukai vector is $v$.
We assume that $M$ is computed with respect to a polarization of $X$
such that $M$ is compact and non-empty, so that $M$ is again a K3.
The integer $n$ is the one defined in~\ref{thm:Muk3}, and $\phi$ is
any of the correspondences defined in~\ref{subsec:MukIsom}. 
\vspace{2mm}

\noindent
{\bf Theorem~\ref{thm:mainthm1}.}  
{\em 
Let $u\in \tH(X, \Z)$ be such that $(u.v) = 1 \bmod n$.  Then
$[\phi(u)]\in\Br(M)$ is the obstruction to the existence of a universal
sheaf on $X\times M$, as defined in~\ref{subsec:quSheaf}.
}

\begin{remark}
Note that since $\tH(X, \Z)$ is unimodular, an $u$ with $(u.v) =1\bmod
n$ can always be found.
\end{remark}

\begin{proof}
First, assume that the moduli problem is fine, so that $\cP$ and
$\cQ$ are untwisted and $n=1$.  We can therefore consider the
correspondence 
\[ \phi = \fmXM^{\cQ, \cO_M}:\tH(X, \Z) \ra \tH(M, \Z). \]
Since $\phi$ is an isometry, and it maps $v$ to $(0,0,\omega)$, it
follows that the degree 0 part of $\phi(u)$ is precisely $(u.v)$.
Using~(\ref{subsec:chclass}) and the projection formula, we get that
for any locally free sheaf $\cE$ on $M$ we have
\begin{align*}
\fmXM^{\cQ, \cE}(u) & = \phi(u)
. \frac{\ch(\cE)}{\rk(\cE)},
\intertext{and therefore the $H^2(M, \Q)$ component of $\fmXM^{\cQ,
\cE}(u)$ satisfies}
\fmXM^{\cQ, \cE}(u)_2 & = \phi(u)_0 \frac{c_1(\cE)}{\rk(\cE)} + \phi(u)_2 \\
& = (u.v) \frac{c_1(\cE)}{\rk(\cE)} + \mbox{ integral part.}
\end{align*}

Let's move on now to the case when the moduli problem is not fine.  We
want to proceed by deforming $X$ until the problem becomes fine, and
this can be done by the argument in~\cite[pp. 385-386]{Muk}.  More
precisely, we can find a smooth family $\cX\ra T_0$ over a small
analytic disk $T_0$, with the following properties:
\begin{enumerate}
\item There is a distinguished point $1\in T_0$ such that $\cX_1$ is
isomorphic to $X$; we'll identify $\cX_1$ with $X$ from here on.
\item The restriction homomorphisms $H^i(\cX, \Z) \ra H^i(\cX_t, \Z)$
are isomorphisms for all $t\in T_0$, and so the cohomology groups of
all the fibers are naturally identified.
\item The Mukai vector $v$ from $\cX_1$ is algebraic in the
Mukai lattice of each fiber $\cX_t$.
\item The polarization of $\cX_1$ is algebraic and ample in each fiber 
$\cX_t$, and therefore all the fibers are naturally polarized.
\item For each fiber $\cX_t$, the moduli space $M(\cX_t, v)$ is
compact and non-empty when computed with respect to this natural
polarization, so it is a K3 surface.  The family of relative moduli
spaces, $\cM\ra T_0$, is smooth over $T_0$.
\item There is a distinguished point $0\in T_0$ such that $M(\cX_0, v)$
is a {\em fine} moduli space of sheaves on $\cX_0$.
\item There exists a twisting $\alpha$ on $\cM$, and a 
$\pi_\cM^*\alpha^{-1}$-twisted sheaf $\cP$ on $\cX\times_{T_0} \cM$,
which is flat over $\cM$, and which restricts to a twisted universal
sheaf on $\cX_t\times\cM_t$ for each $t\in T_0$.
\end{enumerate}

On $\cM$ there exists an $\alpha^{-1}$-twisted locally free sheaf
$\cE$: one can take, for example,
\[ \cE = \pi_{\cM,*} (\pi_\cX^*\cO_\cX(n) \otimes \cP), \]
for a sufficiently high multiple $\cO_\cX(n)$ of a relative
polarization of $\cX/T_0$ (use the flatness of $\cP$ over $\cM$).
Using $\cE$, we can define a global correspondence 
\[ \phi^{\cQ, \cE}:\tH(\cX, \Q) \ra \tH(\cM, \Q), \]
which restricts to the correspondence
\[\phi_{\cX_t \ra \cM_t}^{\cQ_t, \cE_t}:\tH(\cX_t, \Q) \ra \tH(\cM_t,
\Q) \]
for each $t\in T_0$.  Note that since the groups in question are
discrete, these correspondences are necessarily constant as $t$ varies
in $T_0$.

We are now in the situation of Theorem~\ref{thm:localpha}: $\cE$ is a
locally free $\alpha^{-1}$-twisted sheaf on $\cM$, and
$\alpha|_{\cM_0} = 0$ because at $t=0$ the moduli problem is fine.
Therefore, under the corresponding identifications, $\alpha|_{\cM_1}$
is the image of $c_1(\cE_0)/\rk(\cE_0)$ in $\Br(\cM_1)$, where $\cE_0$
is a gluing of $\cE|{\cM_0}$ to a locally free untwisted sheaf.  By
the calculation in the beginning of the proof, we have
\[ \phi_{\cX_0\ra \cM_0}^{\cQ_0, \cE_0}(u)_2 = 
(u.v)\frac{c_1(\cE_0)}{\rk(\cE_0)} + \mbox{ integral part.} \]
But since the correspondences $\phi_{\cX_t \ra \cM_t}^{\cQ_t, \cE_t}$
are constant as a function of $t$, we also have
\[ \phi_{\cX_1 \ra \cM_1}^{\cQ_1, \cE_1}(u)_2 = 
(u.v)\frac{c_1(\cE_0)}{\rk(\cE_0)} + \mbox{ integral part.} \]
Mapping to $\Br(\cM_1)$, we get
\[ [\phi(u)] = [\phi_{\cX_1 \ra \cM_1}^{\cQ_1, \cE_1}(u)_2] = 
\left [ (u.v)\frac{c_1(\cE_0)}{\rk(\cE_0)} + 
\mbox{ integral part}\right ] = (u.v) \alpha, \]
and since $\alpha$ is $n$-torsion, the assumption that $(u.v) = 1 \bmod n$
implies that 
\[ [\phi(u)] = \alpha, \]
which is what we wanted.

To finish the proof, we only need to prove that the choices that we
have made in the above proof do not matter: the actual proof shows
that the choice of $u$ is irrelevant, and the fact that the choice of
$\cE$ is irrelevant (for example, one may choose an $\cE$ on $M$ that
does not extend to the full family $\cM$)
is~(\ref{subsec:phitrindep}).
\end{proof}

\subsection{}
Now let's move on to the proof of Theorem~\ref{thm:mainthm2}.  To
prove that $\alpha=[\phi(u)]$ is in $\Ker \phi^\chk$,
using~(\ref{subsec:defproj}) we need to show that
\[ (\phi(u).\phi(t)) \in \Z~\mbox{ for all } t\in T_X. \]
But since $\phi$ is an isometry, the above is equivalent to 
\[ (u.t) \in\Z\mbox{ for all } t\in T_X, \]
which is obvious.

Let $\lambda\in T_X$ be such that $v-\lambda$ is divisible by $n$ in
$\tH(X, \Z)$ (Theorem~\ref{thm:Muk3}).  Then $(u.\lambda) = (u.v) =
1\bmod n$, and hence $(\phi(u).\phi(\lambda)) = 1 \bmod n$.  But
$\phi(\lambda)$ is divisible by $n$ in $T_M$, $\phi(\lambda) =
n\lambda'$, so we conclude that
\[ (\phi(u). \lambda') = \frac{1}{n} + \mbox{ integer}. \]
This implies that $\alpha = [\phi(u)]\in\Br(M)$ has order at least
$n$.  But by Theorem~\ref{thm:Muk3} $\alpha$ is $n$-torsion, and the
kernel of $\phi^\chk$ is cyclic of order $n$, so we conclude that
$\alpha$ generates $\Ker \phi^\chk$ which is part (1) of
Theorem~\ref{thm:mainthm2}.

This implies at once the equality $\Ker \alpha = \phi(T_X)$ (and not
just $\phi(T_X)\subseteq \Ker\alpha$).  Since $\phi$ is a Hodge
isometry $\tH(X, \Q) \ra \tH(M, \Q)$, it follows that $\phi$ restricts
to a Hodge isometry $T_X\iso \phi(T_X) = \Ker\alpha$, which is part
(2).

Finally, let $\cP$ be a $\pi_M^*\alpha^{-1}$-twisted universal sheaf
on $X\times M$.  To show that 
\[ \FMMX^\cP:\D(M, \alpha) \ra \D(X) \]
is an equivalence of categories, we need to verify the conditions in
Theorem~\ref{thm:equiv}.  For $m\in M$, let $\cP_m$ be the stable sheaf
on $X$ that corresponds to $m$.  The condition
\[ \Hom_X(\cP_m, \cP_m) = \C, \]
follows from the fact that a stable sheaf is simple.  If $m_1\neq
m_2$, $\cP_{m_1}\not\iso \cP_{m_2}$, and they are both stable, so 
\[ \Hom_X(\cP_{m_1}, \cP_{m_2}) = 0. \]
By Serre duality
\[ \Ext^2_X(\cP_{m_1}, \cP_{m_2})  = 0, \]
and since
\[ \chi(\cP_{m_1}, \cP_{m_2}) = -(v,v) = 0, \]
it follows that 
\[ \Ext^1_X(\cP_{m_1}, \cP_{m_2})  = 0. \]
Since $\cP_m$ is a sheaf on $X$ for all $m\in M$ (and not a complex of
sheaves), and $\omega_X = \cO_X$, the remaining conditions of
Theorem~\ref{thm:equiv} are satisfied, and therefore $\FMMX^\cP$ is an
equivalence of categories
\[ \D(M, \alpha) \iso \D(X). \]

\subsection{}
We conclude with a discussion of what the obstacles are to proving
Conjecture~\ref{conj:mainconj}.  Given $X, Y$ and $\alpha, \beta$ as
in the statement of the conjecture, assume that $\Ker\alpha$ is Hodge
isometric to $\Ker\beta$.  A na\"\i ve approach would be to try to
find a third K3 $Z$, with $T_Z \iso \Ker\alpha\iso \Ker\beta$, and to
try to realize $X$ and $Y$ as moduli spaces of stable sheaves on $Z$,
with obstructions $\alpha$ and $\beta$, respectively.  This would
yield $\D(Z) \iso \D(X, \alpha)$ and $\D(Z) \iso \D(Y,
\beta)$, which by transitivity would give the desired equivalence
$\D(X, \alpha) \iso
\D(Y, \beta)$.  However, this approach is soon seen to be too
simplistic: while any sublattice $T$ of $T_X$ such that $T_X/T$ is
cyclic can occur as $\Ker\alpha$ for some $\alpha$, not all such $T$
can be primitively embedded in $\cL_\Kt$.  In other words not all
$\Ker\alpha$ can occur as the transcendental lattice $T_Z$ of some K3
$Z$.  In a certain sense, considering only moduli spaces of {\em
untwisted} sheaves is too restrictive.

The solution to this seems to be the following: one would need to
define a notion of stability for twisted sheaves, and consider
moduli spaces of stable sheaves with arbitrary twisting, instead of
just untwisted ones.  In this case, the universal sheaf would be
twisted in two directions, one from the space $X$ where the stable
sheaves live, and the other from the moduli space $M$, as the
obstruction to the existence of a universal sheaf.  With the extra
flexibility available, one could hope to fully bypass the extra space
$Z$, and to be able to view $Y$ as a moduli space of stable
$\alpha$-twisted sheaves on $X$, with $\beta$ being the obstruction to
the existence of a universal sheaf on $X\times Y$.  This would fully
reproduce Orlov's picture from the untwisted situation.

\subsection{}
The first step thus seems to be defining an appropriate notion of
stability for twisted sheaves.  There is a natural way of doing this:
fix a space $X$ of dimension $d$, polarized by means of a very ample
line bundle of first Chern class $H$, and fix an $\alpha^{-1}$-twisted
locally free sheaf $\cE$.  Then we'll say that an $\alpha$-twisted
sheaf $\cF$ is slope-stable if and only if for every non-trivial
subsheaf $\cG\subset\cF$ we have
\[ \mu_\cE(\cG) < \mu_\cE(\cF), \]
where 
\[ \mu_\cE(\cF) = \frac{\deg_\cE(\cF)}{\rk(\cF)}, \]
$\deg_\cE(\cF)$ being defined as $c_{1, \cE}(\cF).H^{d-1}$, with the
Chern class defined by means of $\cE$ as in~\ref{subsec:chclass}.  (A
similar definition can be given for Maruyama stability of general
sheaves.)

Note that although this definition closely mimics the corresponding
one for untwisted sheaves, we are forced to use the extra data of the
locally free sheaf $\cE$.  In the untwisted context, this has recently
been studied by Yoshioka~(\cite{Yos}); the fact that he was led to the
same definition, coming from a different problem, suggests that this
should indeed be the correct way of approaching stability of twisted
sheaves.

The next step is then to study properties of stable twisted
sheaves, construct moduli spaces, and redo the work of Mukai and
Orlov in the twisted context.  We leave this for a future paper.

\bigskip \noindent
\small\textsc{Department of Mathematics and Statistics, \\
University of Massachusetts, \\
Amherst, MA 01003-4515, USA} \\
{\em e-mail: }{\tt andreic@math.umass.edu}

\end{document}